\documentclass[12pt,reqno]{amsart}

\usepackage{amsmath}\allowdisplaybreaks  
\usepackage{amssymb}        
\usepackage{amsfonts}       
\usepackage{amsthm}         
\usepackage{mathrsfs}       
\usepackage{cases}          
\usepackage{siunitx}        

\usepackage{array}
\usepackage{booktabs}       
\usepackage{longtable}      
\usepackage{multirow}
\usepackage{makecell}       
\usepackage{tabularx}       
\usepackage{nicematrix}     

\usepackage{tikz}   
\usepackage{float} 
\usepackage{xcolor}           
\usepackage{graphicx}         

\usepackage[linesnumbered, ruled, vlined]{algorithm2e}  

\usepackage[colorlinks=true, linkcolor=blue, citecolor=blue, urlcolor=blue]{hyperref} 
\usepackage[nameinlink]{cleveref}	

\usepackage[left=2cm, right=2cm, top=2cm, bottom=2cm]{geometry}

\usepackage{ulem}           

\makeatletter
\newcommand\figcaption{\def\@captype{figure}\caption} 
\newcommand\tabcaption{\def\@captype{table}\caption}
\makeatother

\newtheorem{theorem}{Theorem}[section]          
\newtheorem{example}[theorem]{Example}          

\begin{document}
    \title[RINN for PDEs]{Rank Inspired Neural Network for solving linear partial differential equations}  
    \author[W. Peng, Y. Huang and N. Yi]{Wentao Peng$^1$ and Yunqing Huang$^{1, 2, *}$ and Nianyu Yi$^{1, 3}$ }
    \address{$^1$ School of Mathematics and Computational Science, Xiangtan University, Xiangtan 411105, P.R.China }
\address{$^2$ National Center for Applied Mathematics in Hunan, Xiangtan 411105, P.R.China}
\address{$^3$ Hunan Key Laboratory for Computation and Simulation in Science and Engineering, Xiangtan University, Xiangtan 411105, P.R.China}

\email{aaronpeng98@163.com (W. Peng);\ huangyq@xtu.edu.cn (Y. Huang);\ yinianyu@xtu.edu.cn (N. Yi)}

    
    \begin{abstract}

    This paper proposes a rank inspired neural network (RINN) to tackle the initialization sensitivity issue of physics informed extreme learning machines (PIELM) when numerically solving partial differential equations (PDEs). Unlike PIELM which randomly initializes the parameters of its hidden layers, RINN incorporates a preconditioning stage. In this stage, covariance-driven regularization is employed to optimize the orthogonality of the basis functions generated by the last hidden layer. The key innovation lies in minimizing the off-diagonal elements of the covariance matrix derived from the hidden-layer output. By doing so, pairwise orthogonality constraints across collocation points are enforced which effectively enhances both the numerical stability and the approximation ability of the optimized function space.
    The RINN algorithm unfolds in two sequential stages. First, it conducts a non-linear optimization process to orthogonalize the basis functions. Subsequently, it solves the PDE constraints using linear least-squares method. Extensive numerical experiments demonstrate that RINN significantly reduces performance variability due to parameter initialization compared to PIELM. Incorporating an early stopping mechanism based on PDE loss further improves stability, ensuring consistently high accuracy across diverse initialization settings.
    \end{abstract}
\thanks{$^*$ Corresponding author.}
    
    \maketitle

\section{Introduction}
    Partial differential equations (PDEs) play a fundamental role in modeling physical phenomena and engineering systems. Classical grid‑based numerical methods, including the finite element method \cite{rao2010finite}, finite difference method \cite{leveque2007finite} and finite volume method \cite{versteeg2007introduction}, have been extensively validated and widely applied. As an alternative, deep neural networks have in recent years emerged as a mesh‑free, end‑to‑end framework for PDE solution due to their flexible representational capacity. Lagaris et al.\ \cite{lagaris1998artificial,lagaris2000neural} were the first to apply multilayer perceptrons as global collocation approximators for PDEs, drawing on Cybenko’s universal approximation theorem \cite{cybenko1989approximation}. Subsequently, weak‑form approaches such as the Deep Ritz Method \cite{yu2018deep}, Deep Nitsche Method \cite{ming2021deep} and Weak Adversarial Networks \cite{zang2020weak} recast the variational formulation of a PDE into a loss functional and approximate weak solutions via automatic differentiation or adversarial training, while strong‑form approaches exemplified by Physics‑Informed Neural Networks (PINNs) \cite{raissi2019physics} and the Deep Galerkin Method (DGM) \cite{sirignano2018dgm} enforce the PDE residual, boundary and initial conditions pointwise, offering straightforward implementation but often suffering from spectral bias \cite{rahaman2019spectral} and instability in high‑frequency regimes.  
    
    Since their introduction, PINNs have been extended in several directions. vPINNs \cite{kharazmi2019variational} and hp‑vPINNs \cite{kharazmi2021hp} enhance stability through variational test spaces. Adaptive sampling \cite{wu2023comprehensive} and domain decomposition schemes \cite{jagtap2020conservative,jagtap2020extended} accelerate convergence in regions with localized features. Meanwhile, transfer‑learning extensions \cite{liu2023adaptive} broaden the framework’s applicability. Qian et al.\ \cite{qian2023physics} analyzed the application of PINNs to three second-order hyperbolic partial differential equations (PDEs). By introducing a modified loss function, they enhanced PINNs' ability to approximate dynamic problems and subsequently proved the method's generalization error bound. This bound demonstrates that the approximation errors for the solution field, its time derivative, and spatial gradient are jointly governed by the training loss magnitude and the number of orthogonal sampling points employed. However, PINNs still suffer from premature convergence to local minima and exhibit difficulty in learning high-frequency oscillatory features.
    
    An alternative mesh-free method is provided by the extreme learning machine (ELM) family \cite{huang2006extreme,huang2006universal}, where only the output weights are trained and the hidden layer parameters are randomly fixed. Dwivedi et al.\ \cite{dwivedi2020physics} proposed the Physics Informed Extreme Learning Machine (PIELM) which combines ELM efficiency with the physics-based loss formulation of PINNs. By randomly initializing the hidden-layer parameters and solving a physics-constrained least-squares problem for the output weights, PIELM achieves fast and accurate solutions for linear PDEs without iterative optimization. Initially limited to linear PDEs, PIELM was extended to handle nonlinear problems by Dong et al.\ \cite{dong2021local} using specialized solvers that directly optimize the full nonlinear system via perturbation strategies and Newton-least squares methods. Shang et al.\ \cite{shang2022deep} combined local randomization neural networks with the Petrov–Galerkin framework. Sun et al. \cite{sun2024local} employed randomized neural networks as trial functions within a discontinuous Galerkin framework (DG-RNN) to solve nonlinear PDEs, enhancing shock-capturing capability and convergence order. Wang and Dong \cite{wang2024extreme} proposed randomized neural networks for high-dimensional PDEs, introducing the A-TFC algorithm to reduce boundary condition complexity from exponential to linear scaling.

    Although PIELM and its variants offer simplicity and rapid training, solution accuracy depends critically on hidden-layer initialization. Dong et al.\ \cite{dong2022computing} studied the influence of initialization on PIELM and proposed an optimization framework based on residuals. This framework uses differential evolution to identify the hidden layer hyperparameters that are close to the optimal values. Calabro et al.\ \cite{calabro2021extreme} comprehensively investigated weight initialization strategies for PIELM, establishing theoretically grounded parameter ranges that enable robust solutions to one-dimensional elliptic equations with sharp boundary/internal layers.

   Inspired by the orthogonality principles of spectral methods, a Rank‑Inspired Neural Network (RINN) framework is proposed, significantly enhancing the stability of PIELM by enforcing near‑orthogonality among hidden‑layer outputs prior to coefficient solving.   Our primary contributions are:
    \begin{itemize}
      \item \textbf{Covariance‐driven basis orthogonalization.} A pre‑training stage is introduced to minimize the off‑diagonal elements of the covariance matrix of hidden‑layer outputs, thereby reducing the correlation between neural basis functions and alleviating the initialization sensitivity of ELM.
      \item \textbf{Algorithm early-stopping mechanism.} In each iteration, hidden‑layer parameters are optimized through covariance‑driven basis function orthogonalization. Upon completion, a system of linear equations associated with the PDE is constructed to solve for the output‑layer parameters. The PDE residual loss is then computed based on these parameters and training is stopped according to the observed changes in this loss.
    \end{itemize}
    
    The remainder of this paper is organized as follows. \Cref{section:2} reviews the PIELM algorithm and highlights its sensitivity to initialization. Building on this, a covariance‑driven basis function orthogonalization method is presented, leading to the Rank‑Inspired Neural Network (RINN) algorithm. A series of numerical experiments then verify RINN’s effectiveness. \Cref{section:3} demonstrates that excessive training can negatively affect the performance of RINN. Therefore, RINN with early stopping (RINN-es) was proposed and its effectiveness was verified in subsequent experiments. Finally, \Cref{section:4} concludes the paper and outlines directions for future research.
    
\section{Rank Inspired Neural Network}
\label{section:2}
    In this section, we consider the following linear partial differential equation, 
    \begin{equation}
        \left\{
        \begin{aligned}
            \mathscr{A}u(\boldsymbol{x},t) &= f(\boldsymbol{x},t), && \text{in } \Omega \times I, \\
            \mathscr{B}u(\boldsymbol{x},t) &= g(\boldsymbol{x},t), && \text{on } \partial \Omega \times I, \\
            u(\boldsymbol{x},0) &= h(\boldsymbol{x}),   && \text{in } \Omega \times \{0\},
        \end{aligned}
        \right.
        \label{eq1}
    \end{equation}
    Where $\Omega \subset \mathbb{R}^d$ denotes the spatial domain with boundary $\partial \Omega$, and $I = [0, T]$ ($T > 0$) represents the time interval. The variables $\boldsymbol{x} = (x_1, \dots, x_d) \in \Omega$ and $t \in I$ refer to the spatial and temporal coordinates, respectively. $\mathscr{A}$ and $\mathscr{B}$ are linear differential operators. The function $u(\boldsymbol{x}, t)$ is the unknown solution, while $f(\boldsymbol{x}, t)$ is a prescribed source term and $g(\boldsymbol{x}, t)$ specifies the boundary condition.
    
\subsection{Physics Informed Extreme Learning Machine}
	In this subsection, we introduce the PIELM for the proposed equation \eqref{eq1}. 
	
\begin{figure}[htbp]
    	\centering
    	\includegraphics[width=0.5\linewidth]{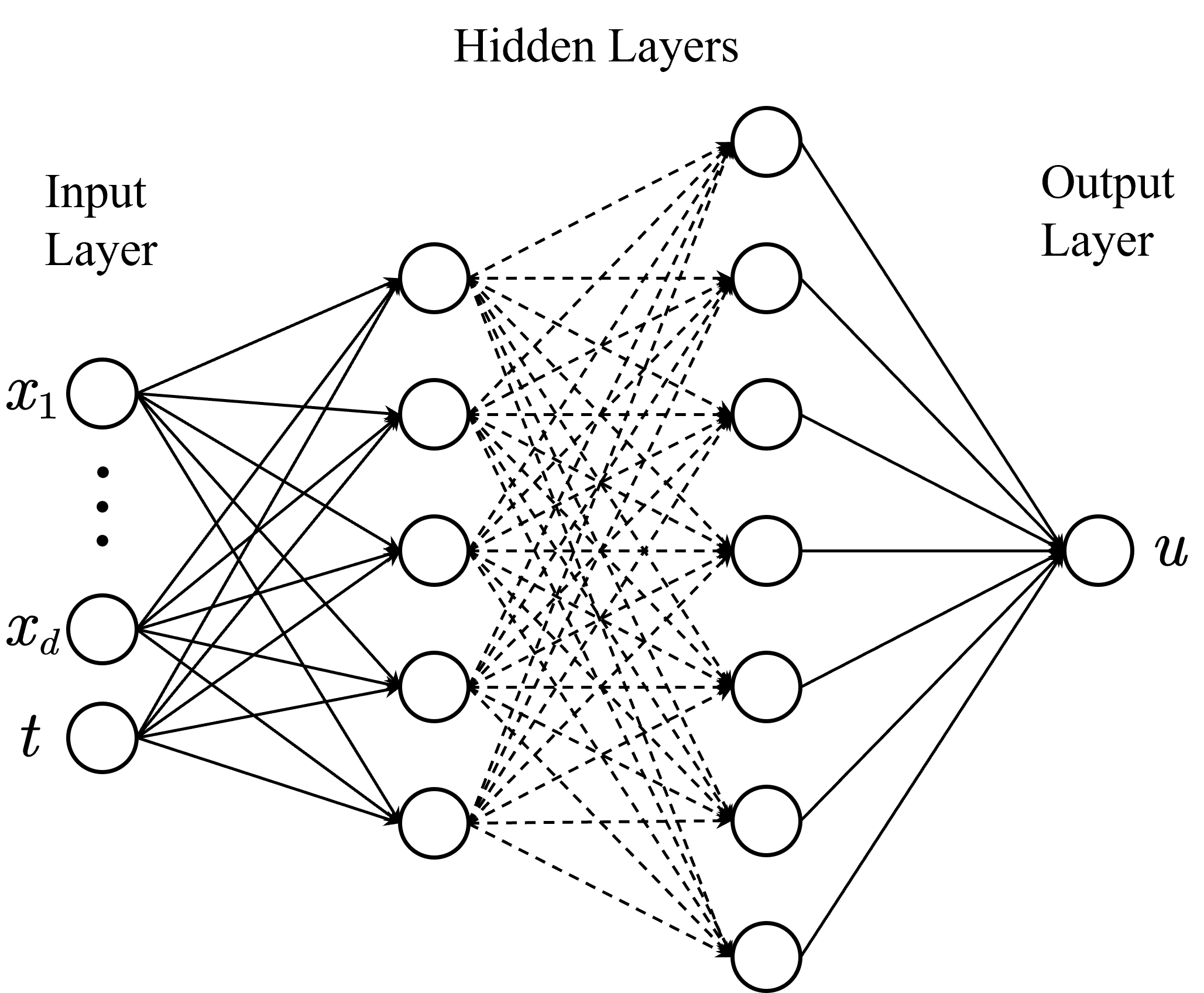}
        \caption{Structure of ELM.}
        \label{fig:net}
    \end{figure}

    The solution $u(\boldsymbol{x})$ of the equation \eqref{eq1} is approximated using an extreme learning machine (ELM) with at least one hidden layer ($L \geq 1$), where only the output layer is trained. A standard network topology of ELM is shown in \Cref{fig:net}.

    The input layer contains $d+1$ neurons which correspond to spatial coordinates $\boldsymbol{x} = (x_1,\dots,x_d)$ and temporal coordinate $t$. Each hidden layer, indexed by $\ell \in \{1,\dots,L\}$, receives the output from the previous layer and processes it using a linear transformation followed by a nonlinear activation function. This operation is defined as:
    
    \begin{equation}
    \boldsymbol{z}_{\ell} = \sigma\left( \boldsymbol{W}^{(\ell)}\boldsymbol{z}_{\ell-1} + \boldsymbol{b}^{(\ell)} \right), \qquad \boldsymbol{W}^{(\ell)} \in \mathbb{R}^{N_\ell \times N_{\ell-1}}, \ \boldsymbol{b}^{(\ell)} \in \mathbb{R}^{N_\ell},
    \end{equation}
     where $\sigma(\cdot)$ denotes the activation function, $\boldsymbol{W}^{(\ell)}$ represents the weight matrix, and $\boldsymbol{b}^{(\ell)}$ is the bias vector. The collection of all these parameters, $\boldsymbol{\theta} = \{\boldsymbol{W}^{(1)},\boldsymbol{b}^{(1)},\dots,\boldsymbol{W}^{(L)},\boldsymbol{b}^{(L)}\}$, remains fixed during training. The total number of parameters is $N_{\boldsymbol{\theta}} = \sum_{\ell=1}^L N_\ell(N_{\ell-1} + 1)$ which depends on the network structure $[N_0 = d+1, N_1,\dots,N_L, N_{L+1} = 1]$. 
    
    Define the output of each neuron in the final hidden layer as a kind of neural basis function. The $j$-th basis function is given by:
    
    \begin{equation}
    \phi_j(\boldsymbol{\theta},\boldsymbol{x}) = \left[ \boldsymbol{z}_L \circ \cdots \circ \boldsymbol{z}_1 \right]_j (\boldsymbol{x}).
    \label{eq:basis}
    \end{equation}
    
    The complete output of the last hidden layer can be written as a vector of these functions:
    $$
    \boldsymbol{\Phi}(\boldsymbol{\theta},\boldsymbol{x}) = [\phi_1(\boldsymbol{\theta},\boldsymbol{x}),\dots,\phi_{N_L}(\boldsymbol{\theta},\boldsymbol{x})].
    $$
    To distinguish the trainable components, we enforce a linear activation function in the output
    layer and denote the weights of the output layer as $\boldsymbol{\beta} = (\beta_1,\dots,\beta_{N_L})^\top$, then the final output of the ELM is a linear combination of the basis functions $\phi_j(\boldsymbol{\theta},\boldsymbol{x})$ weighted by $\beta_j$, for $j=1,\dots,N_L$:    
    
    \begin{equation}
    u(\boldsymbol{x}) = \sum_{j=1}^{N_L} \beta_j \phi_j(\boldsymbol{\theta},\boldsymbol{x}) = \boldsymbol{\Phi}(\boldsymbol{\theta},\boldsymbol{x})\boldsymbol{\beta}.
    \label{ELM}
    \end{equation}
	 
	The network input consists of a set of collocation points sampled from the spatio-temporal domain. We define three distinct sets of points:

    \begin{itemize}
        \item Residual points: $\mathbf{X}_{\mathrm{res}} \in \mathbb{R}^{K_{\mathrm{res}} \times d}$ in $\Omega \times I$.
        \item Boundary condition points: $\mathbf{X}_{\mathrm{bcs}} \in \mathbb{R}^{K_{\mathrm{bcs}} \times d}$ on $\partial \Omega \times I$. 
        \item Initial condition points: $\mathbf{X}_{\mathrm{ics}} \in \mathbb{R}^{K_{\mathrm{ics}} \times d}$ on $\Omega \times \{0\}$.
    \end{itemize}
    All of these points together form the complete set of collocation points, denoted by $\mathcal{X} = \{\mathbf{X}_{\mathrm{res}}, \mathbf{X}_{\mathrm{bcs}},\allowbreak \mathbf{X}_{\mathrm{ics}}\}$, with a total of $K = K_{\mathrm{res}} + K_{\mathrm{bcs}} + K_{\mathrm{ics}}$ points. These points are assembled into a global input matrix $\mathbf{X} \in \mathbb{R}^{K \times d}$ where each row $x_k$ corresponds to a specific space-time coordinate.
	 
	We apply the ELM \eqref{ELM} with the collocation points $\mathbf{X}_{\mathrm{res}},\mathbf{X}_{\mathrm{bcs}},\mathbf{X}_{\mathrm{ics}} \subset \mathcal{X}$ to discretize equation \eqref{eq1}, leading to the following linear system:

    \begin{equation}
        \left\{
        \begin{aligned}
            \sum_{i=1}^N [\mathscr{A}\phi_i(\boldsymbol{\theta}, \mathbf{X}_{\mathrm{res}})] \beta_i  
	 		&= f(\mathbf{X}_{\mathrm{res}}) ,\quad &&\mathbf{X}_{\mathrm{res}} \in \Omega \times I,  \\
            \sum_{i=1}^N [\mathscr{B}\phi_i(\boldsymbol{\theta}, \mathbf{X}_{\mathrm{bcs}})] \beta_i  
	 		&= g(\mathbf{X}_{\mathrm{bcs}}) ,\quad &&\mathbf{X}_{\mathrm{bcs}} \in \partial \Omega \times I, \\
            \sum_{i=1}^N [\phi_i(\boldsymbol{\theta}, \mathbf{X}_{\mathrm{ics}})] \beta_i  
	 		&= h(\mathbf{X}_{\mathrm{ics}}) ,\quad &&\mathbf{X}_{\mathrm{ics}} \in \Omega \times \{0\}.
        \end{aligned}
        \right.
        \label{eq3}
    \end{equation}
This results in a system of $K$ algebraic equations involving the trainable parameters $\boldsymbol{\beta}$. For a given $\boldsymbol{\theta}$, the coefficients $\mathscr{A}\phi_i(\boldsymbol{\theta}, \mathbf{X}_{\mathrm{res}}),\mathscr{B}\phi_i(\boldsymbol{\theta}, \mathbf{X}_{\mathrm{bcs}}),\phi_i(\boldsymbol{\theta}, \mathbf{X}_{\mathrm{ics}})$ can be evaluated using forward propagation and automatic differentiation.
	
Since system \eqref{eq3} is linear in $\boldsymbol{\beta}$, it can be obtained by solving a least‑squares problem. For simplicity, the system is then expressed in matrix form as:

	\begin{equation}
	 	\boldsymbol{H}(\boldsymbol{\theta},\mathbf{X}) \boldsymbol{\beta} = \boldsymbol{S},\qquad 
	 	\boldsymbol{H}(\boldsymbol{\theta},\mathbf{X}) = 
	 	\begin{bmatrix}
	 	\mathscr{A}\boldsymbol{\Phi}(\boldsymbol{\theta}, \mathbf{X}_{\mathrm{res}})	\\
	 	\mathscr{B}\boldsymbol{\Phi}(\boldsymbol{\theta}, \mathbf{X}_{\mathrm{bcs}})	\\
	 	\boldsymbol{\Phi}(\boldsymbol{\theta}, \mathbf{X}_{\mathrm{ics}})	
	 	\end{bmatrix},
	 	\quad 
	 	\boldsymbol{S} = 
	 	\begin{bmatrix}
	 		f(\mathbf{X}_{\mathrm{res}})	\\
	 		g(\mathbf{X}_{\mathrm{bcs}})	\\
	 		h(\mathbf{X}_{\mathrm{ics}})
	 	\end{bmatrix}.
	 	\label{eq:matrix}
	\end{equation}
For any fixed $\boldsymbol{\theta}$, the least squares solution for the linear system is given by,
\begin{equation}
	 	\boldsymbol{\beta} = \left[\boldsymbol{H}(\boldsymbol{\theta},\mathbf{X})\right]^\dagger  \boldsymbol{S},
	\end{equation}
where $\boldsymbol{H}^\dagger $ denotes the Moore–Penrose pseudo-inverse of $\boldsymbol{H}$.
	 
Since the hidden layer parameters $\boldsymbol{\theta}$ in an extreme learning machine (ELM) are fixed after initialization, the performance of PIELM is inherently sensitive to the initialization. To show this sensitivity, we conducted an experiment on the two-dimensional Poisson equation with Dirichlet boundary conditions:
    	$$
        \left\{
		\begin{aligned}
		 		-\Delta u &= f, &&\text{in } \Omega = [-1,1]^2,\\
		 		u&=g,&&\text{on } \partial\Omega.
		\end{aligned}
        \right.
		$$
        We take $\sin(2\pi x)\sin(2\pi y)$ as the solution of this problem and employ a network structure $N_{\text{nn}}=[2,512, 1024,1]$. The network weights are initialized from three different uniform distributions: $\mathcal{U}(-0.5,0.5)$, $\mathcal{U}(-1,1)$, and $\mathcal{U}(-1.5,1.5)$. To construct the linear system, we use $K_{\mathrm{res}}=2048$ residual points and $K_{\mathrm{bcs}}=4096$ boundary condition points.

        \begin{table}[htbp]
            \caption{Comparison of PIELM performance under different weight initializations.}
            \label{table:initial sensitity}
            \begin{tabular*}{0.75\textwidth}{@{\extracolsep{\fill}}lccc}
                \toprule
                Initialization    & $\mathcal{U}(-0.5,0.5)$  & $\mathcal{U}(-1,1)$ & $\mathcal{U}(-1.5,1.5)$ \\ 
                \midrule
                $E_{L^2}$  & $1.84\times 10^{-9}$ & $2.86\times 10^{-3}$ & $6.31\times 10^{-1}$  \\ 
                $E_{L^1}$ & $7.30\times 10^{-10}$ & $6.40\times 10^{-4}$ & $2.23\times 10^{-1}$  \\
                \bottomrule
            \end{tabular*}
        \end{table}

        We calculate the relative $L^2$ error ($E_{L^2}$) and mean absolute error ($E_{L^1}$) as follows,
        \begin{equation}
            E_{L^2} = \tfrac{\sqrt{\sum\limits_{i=1}^N \left| u_{\text{true}}(\mathbf{x}_i) - u_{\text{pred}}(\mathbf{x}_i) \right|^2}}{\sqrt{\sum\limits_{i=1}^N \left| u_{\text{true}}(\mathbf{x}_i) \right|^2}}, 
            \qquad
            E_{L^1} = \frac{1}{N} \sum\limits_{i=1}^N \left| u_{\text{true}}(\mathbf{x}_i) - u_{\text{pred}}(\mathbf{x}_i) \right|.
            \label{Errors}
        \end{equation}

        As shown in \Cref{table:initial sensitity}, the accuracy of PIELM is highly dependent on the initial distribution of fixed network parameters. Among the tested distributions, the relative $L^2$ errors span across eight orders of magnitude, indicating that it is highly sensitive to the choice of weight initialization.

\subsection{Rank Inspired Neural Network}
    In this section, we present the RINN algorithm for solving the linear PDE \eqref{eq1}. The training process of the RINN algorithm is divided into two stages. In the first stage, the weight parameters of the hidden layer are trained through an optimization process. This optimization process is specifically designed to ensure that the output function of the final hidden layer complies with the orthogonality constraint, thereby enhancing the numerical stability and approximation ability of the function space. In the second stage, the parameters of the hidden layer are frozen and the parameters of the output layer are determined by the least squares method to solve the PDE.
    
    The first stage, referred to as covariance-driven basis orthogonalization, is motivated by classical Galerkin methods. In such traditional numerical methods, orthogonal basis functions are typically utilized to obtain well-conditioned full-rank linear systems. In contrast, the neural basis functions defined in \eqref{eq:basis} may exhibit significant mutual correlations which can amplify the effects of random initialization and hinder convergence.

	Ideally, the basis functions should be pairwise orthogonal in the $L^2(\Omega)$ sense, that is
    $$
    \int_{\Omega} \phi_i(\boldsymbol{\theta}, \boldsymbol{x})\, \phi_j(\boldsymbol{\theta}, \boldsymbol{x})\, d\boldsymbol{x} = 0 ,\quad  i \ne j.
    $$

    In practice, we approximate this condition using the collocation points in the spatial domain, yielding the discrete orthogonality constraint:
    $$
    \sum_k \phi_i(\boldsymbol{\theta}, \boldsymbol{x}_k)\cdot\phi_j(\boldsymbol{\theta}, \boldsymbol{x}_k)=0,\qquad  i\ne j, \ \boldsymbol{x}_k \in \Omega\times I.
    $$
	
	Let $\{\phi_i(\boldsymbol{\theta}, \boldsymbol{x})\}_{i=1}^{N_L}$ denote the set of neural basis functions evaluated at a collection of $K$ collocation points $\{\boldsymbol{x}_j\}_{j=1}^K\ (K\geq N_L)$. We can get the matrix form of these neural basis functions, 
	\begin{equation}\begin{aligned}
		\boldsymbol{\Phi}(\boldsymbol{\theta}, \boldsymbol{x}) = &(\phi_1(\boldsymbol{\theta}, \boldsymbol{x}), \phi_2(\boldsymbol{\theta}, \boldsymbol{x}), \dots, \phi_{N_L}(\boldsymbol{\theta}, \boldsymbol{x})) \\
        =& 
		\begin{bmatrix}
			\phi_1(\boldsymbol{\theta}, \boldsymbol{x}_1) & \phi_2(\boldsymbol{\theta}, \boldsymbol{x}_1) & \cdots & \phi_{N_L}(\boldsymbol{\theta}, \boldsymbol{x}_1)\\
			\phi_1(\boldsymbol{\theta}, \boldsymbol{x}_2) & \phi_2(\boldsymbol{\theta}, \boldsymbol{x}_2) & \cdots & \phi_{N_L}(\boldsymbol{\theta}, \boldsymbol{x}_2)\\
			\vdots      & \vdots      & \ddots & \vdots\\
			\phi_1(\boldsymbol{\theta}, \boldsymbol{x}_K) & \phi_2(\boldsymbol{\theta}, \boldsymbol{x}_K) & \cdots & \phi_{N_L}(\boldsymbol{\theta}, \boldsymbol{x}_K)
		\end{bmatrix},
        \end{aligned}
	\end{equation}
    where the $j$-th row corresponds to the value of all neural basis functions at point $\boldsymbol{x}_j$. The primary objective is to potimize $\boldsymbol{\theta}$ to maximize the rank of the matrix $\boldsymbol{\Phi}(\boldsymbol{\theta}, \boldsymbol{x})$. A maximized rank ensures that the row space of the matrix has the highest possible dimension, which in turn enhances its ability to represent or approximate complex functions effectively.
	
	To quantify the degree of orthogonality among the basis functions, we consider the sample covariance matrix of $\boldsymbol{\Phi} $, given by:
	\begin{equation}
		\mathbf{C} = \frac1{K-1} \boldsymbol{\Phi}^{\top} \boldsymbol{\Phi}.
	\end{equation}
    where
    \begin{equation}
        \mathbf{C}_{ij}=\frac1{K-1} \sum_{k=1}^K \phi_i(\boldsymbol{\theta}, \boldsymbol{x}_k) \phi_j(\boldsymbol{\theta}, \boldsymbol{x}_k),
    \end{equation}
    which serves as a discrete approximation to the $L^2(\Omega)$ inner product between the $i$-th and $j$-th neural basis functions.
    
    If the bases are orthogonal in the $L^2(\Omega)$ sense, the off-diagonal entries of $\mathbf{C}$ should vanish, i.e., $\mathbf{C}_{ij} = 0$ for $i \ne j$. Hence, the degree of orthogonality among the neural basis functions can be directly inferred from the structure of the covariance matrix $\mathbf{C}$.
    To promote orthogonality among the neural basis functions, we minimize the Frobenius norm of the off-diagonal elements of the sample covariance matrix $\mathbf{C}$. The associated loss function is defined as:
	\begin{equation}
		\mathcal{L}_{\rm ortho} = \|\mathbf{C}\odot(\mathbf{J}-\mathbf{E})\|_F = \sqrt{\sum_{i\ne j} |\mathbf{C}_{ij}|^2}, \label{eq:loss}
	\end{equation}
	where $\mathbf{J}\in\mathbb{R}^{N_L\times N_L}$ denotes the all-ones matrix, $\mathbf{E}\in\mathbb{R}^{N_L\times N_L}$ is the identity matrix, and $\odot$ represents the Hadamard (elementwise) product. Here, $\|\cdot\|_F$ denotes the Frobenius norm and $\mathbf{C}_{ij}$ denotes the $(i,j)$-th element of the covariance matrix $\mathbf{C}$.

    Additionally, we hope that each basis function is normalized to the unit norm, that is
    \begin{equation}
        \|\phi_i\|_{L^2(\Omega)}^2 = \int_{\Omega} \phi_i(\boldsymbol{\theta}, \boldsymbol{x})^2 d\boldsymbol{x}=1.
    \end{equation}
    For a covariance matrix, the diagonal elements must equal 1, i.e., $\mathbf{C}_{ii} = 1$ for $i = 1, 2, \dots, N_L$. To measure and punish the degree to which these diagonal elements deviate from 1, we can define the diagonal element loss as:
    \begin{equation}
        \mathcal{L}_{\mathrm{diag}} = \sum_{i=1}^{N_L} \left| \log_{10} \left( \mathbf{C}_{ii}^{2} \right) \right| = \sum_{i=1}^{N_L} \left | \log_{10}(\mathbf{C}_{ii}) \times 2 \right | \, .
    \end{equation}
    Here, the absolute value of $\log_{10} \left( \mathbf{C}_{ii}^{2} \right)$ is taken because no matter whether $\mathbf{C}_{ii}$ is too large or too small, $\log_{10} \left( \mathbf{C}_{ii}^{2} \right)$ will deviate from 0 and the more it deviates from 0, the stronger the penalty is.

    The total loss function is defined as a weighted sum of the orthogonal loss $\mathcal{L}_{\mathrm{ortho}}$ and the diagonal element loss $\mathcal{L}_{\mathrm{diag}}$:
    \begin{equation}
        \mathcal{L}_{\mathrm{total}} = \varepsilon\mathcal{L}_{\mathrm{diag}} + \mathcal{L}_{\mathrm{ortho}}.
    \end{equation}
    where $\varepsilon > 0$ is a tunable parameter controlling the influence of the diagonal element loss.
    
    This orthogonalization goal is optimized with respect to the trainable parameters defining the neural basis functions. Once the basis has been sufficiently decorrelated, the second stage proceeds by solving the resulting linear system in a least-squares sense to determine the output weights.

    To end this section, we summarize the RINN algorithm for PDEs \eqref{eq1} in Algorithm \ref{alg:RINN1}. 
    
    \begin{algorithm}[H]
        \caption{RINN for Linear PDEs}
        \label{alg:RINN1}
        \KwIn{Collocation points $\mathcal{X} = \{\mathbf{X}_{\mathrm{res}}, \mathbf{X}_{\mathrm{bcs}}, \mathbf{X}_{\mathrm{ics}}\}$, maximum iterations $E$, learning rate $\eta$, regularization coefficient $\varepsilon$.}
        \KwOut{Optimized parameters $(\boldsymbol{\theta}, \boldsymbol{\beta})$.}
        
        Initialize network parameters $\boldsymbol{\theta}$ randomly\;
        
        \textbf{Stage 1: Covariance-Driven Basis Orthogonalization}\;
        \For{$epoch := 1$ \KwTo $E$}{
            Compute hidden layer output matrix $\boldsymbol{\Phi} = \boldsymbol{\Phi}_{\boldsymbol{\theta}}(\mathcal{X})$ via forward propagation\;
            
            Compute sample covariance matrix: $\mathbf{C} = \tfrac{1}{K-1}\, \boldsymbol{\Phi}^\top \boldsymbol{\Phi}$\;
            
            Compute diagonal loss: $\mathcal{L}_{\mathrm{diag}} = \sum_{i=1}^{N_L} \left| \log_{10}(\mathbf{C}_{ii}^2) \right|$\;
            
            Compute off-diagonal loss: $\mathcal{L}_{\mathrm{ortho}} = \|\mathbf{C} \odot (\mathbf{J} - \mathbf{E})\|_F$\;
            
            Compute total loss: $\mathcal{L}_{\mathrm{total}} = \varepsilon\,\mathcal{L}_{\mathrm{diag}} + \mathcal{L}_{\mathrm{ortho}}$\;
            
            Update parameters: $\boldsymbol{\theta} = \boldsymbol{\theta} - \eta \, \nabla_{\boldsymbol{\theta}} \mathcal{L}_{\mathrm{total}}$\;
        }
        
        \textbf{Stage 2: Least-Squares Solve}\;
        
        Assemble linear system $\mathbf{H} \boldsymbol{\beta} = \mathbf{S}$ as defined in \eqref{eq:matrix}\;
        
        Solve for weights: $\boldsymbol{\beta} = \mathbf{H}^\dagger \mathbf{S}$\;
    \end{algorithm}

\subsection{Numerical Examples}

This subsection evaluates RINN through three perspectives: (i) performance on PDEs with analytical solutions, (ii) comparative performance against PIELM, and (iii) computational efficiency analysis. All computations are performed on an NVIDIA GeForce RTX 4080 GPU using PyTorch 3.8.20 with double precision arithmetic, ensuring $15$-digit decimal precision throughout. Our experiments encompass elliptic and evolution equations with analytical solutions.

The hidden layers of all neural networks are activated by the $\tanh$ activation function and the weights are initialized according to a uniform distribution as defined below,
\begin{align*}
\text{Uniform distribution: } & \mathcal{W} \sim \mathcal{U}(-a, a).
\end{align*}

\subsubsection{Elliptic Equations}
Here, we consider the Poisson equation with Dirichlet boundary conditions:
\begin{equation}
\left\{
	\begin{aligned}
		- \Delta u(\boldsymbol{x})&=f(\boldsymbol{x}) ,&&\boldsymbol{x}\in \Omega ,\\
		u(\boldsymbol{x})&=g(\boldsymbol{x}),          &&\boldsymbol{x}\in \partial \Omega .
	\end{aligned}
\right.
	\label{eq:Elliptic Equations}
\end{equation}

\begin{example} 
    Consider 1D Poisson equations defined on the interval $[-1, 1]$ with the following two analytical solutions:
    \begin{subequations}\label{eq:1D_all}
        \begin{align}
          u(x) &= \sin(2\pi x)\cos(3\pi x),
            \label{eq:1D 1}\\
          u(x) &= \sin(7.5\pi x).
            \label{eq:1D 2}
        \end{align}
    \end{subequations}
    \label{ep:Poi1D}
\end{example}

We employ a fully connected neural network with architecture $N_{\text{nn}} = [1, 128, 1]$. The number of collocation points is $K_{\mathrm{res}} = 1024$ and $K_{\mathrm{bcs}} = 2$. All network weights are initialized from the uniform distribution $\mathcal{U}(-20, 20)$.
Training is performed using the Adam optimizer with a learning rate of $\eta = 10^{-3}$ for $E = 2000$ epochs. The regularization parameter in the loss function \eqref{eq:loss} is set to $\varepsilon = 0.1$.



\begin{table}[htbp]
    \caption{Comparison of PIELM and RINN on \Cref{ep:Poi1D}.}
    \label{table:1D test}
    \begin{tabular*}{0.75\textwidth}{@{\extracolsep{\fill}}cccc}
        \toprule
        Analytical Solution   & Method         & {$E_{L^2}$}     & {$E_{L^1}$}    \\
        \midrule
        \multirow{2}{*}{\eqref{eq:1D 1}} 
                    & PIELM & $2.19\times10^{-5}$ & $1.02\times10^{-5}$ \\
                    & RINN  & $1.26\times10^{-8}$ & $5.18\times10^{-9}$ \\
        \midrule
        \multirow{2}{*}{\eqref{eq:1D 2}} 
                    & PIELM & $2.19\times10^{-3}$ & $1.41\times10^{-3}$ \\
                    & RINN  & $6.33\times10^{-7}$ & $3.65\times10^{-7}$ \\
        \bottomrule
    \end{tabular*}
\end{table}

\Cref{table:1D test} reports the performance of PIELM and RINN on the two 1D Poisson equations with analytical solutions. In both cases, RINN consistently outperforms PIELM in terms of accuracy. For the first solution \eqref{eq:1D 1}, RINN reduces the relative $L^2$ error from $\mathcal{O}(10^{-5})$ to $\mathcal{O}(10^{-8})$ and the mean absolute error from $\mathcal{O}(10^{-5})$ to $\mathcal{O}(10^{-9})$. The improvement is even more pronounced in the second solution \eqref{eq:1D 2}, where the relative $L^2$ error decreases from $\mathcal{O}(10^{-3})$ under PIELM to $\mathcal{O}(10^{-7})$ with RINN. 

\begin{figure}[htbp]
	\centering
	\includegraphics[width=0.85\linewidth]{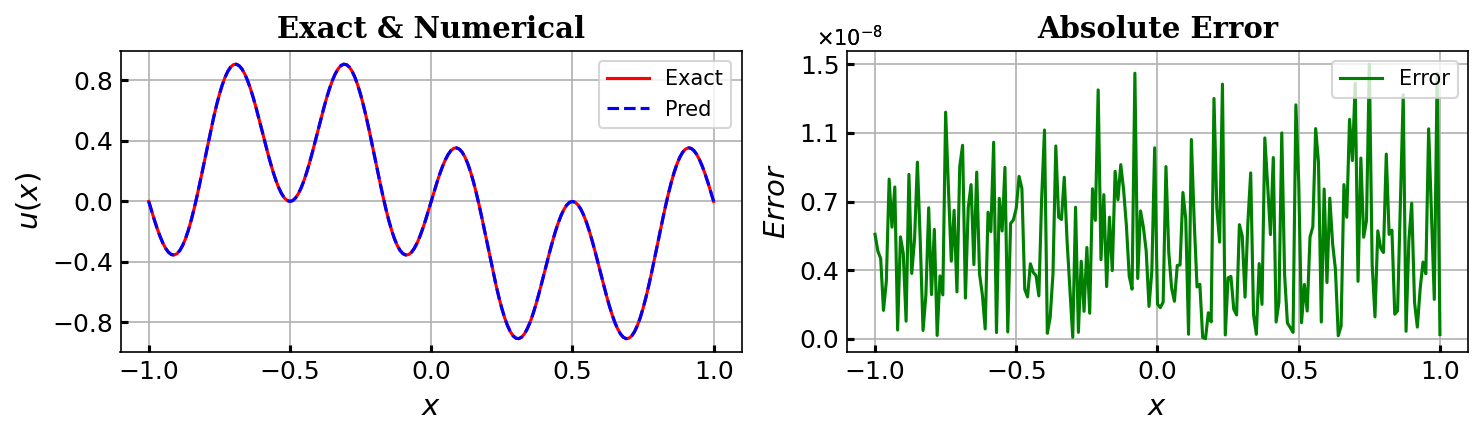}
	
	\includegraphics[width=0.85\linewidth]{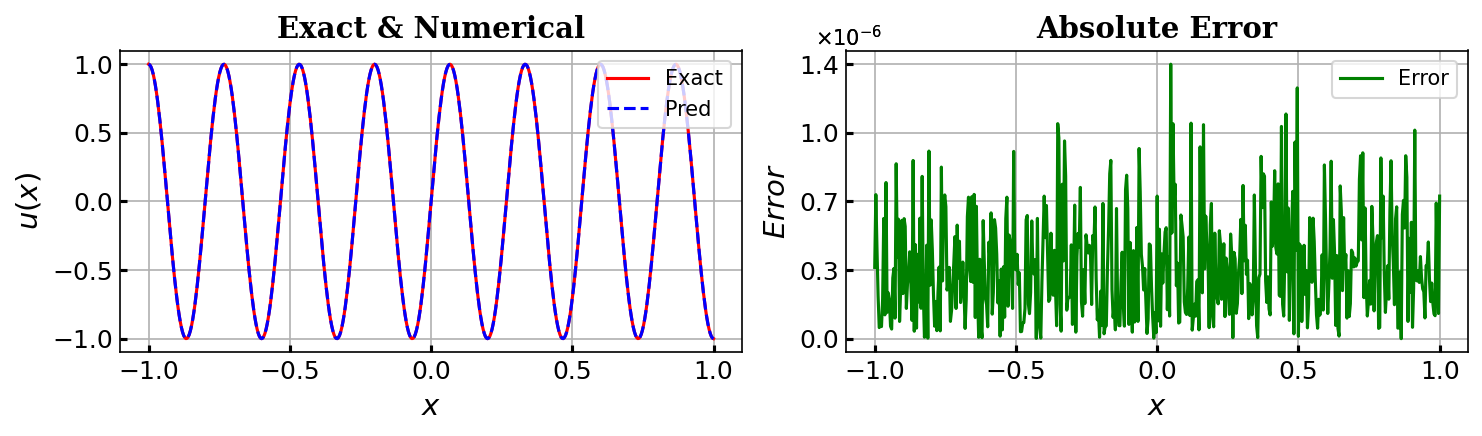}
	\caption{RINN results for \Cref{ep:Poi1D}. Left column shows exact and predicted solution; right column shows pointwise absolute errors. Rows from top to bottom correspond to the solution \eqref{eq:1D 1} and the solution \eqref{eq:1D 2}.}
	\label{fig 1D}
\end{figure}

\Cref{fig 1D} displays the numerical results obtained by RINN for two 1D Poisson equations with distinct exact solutions. In each panel, the left column compares the exact and RINN solutions, while the right column shows the corresponding pointwise absolute errors. The top and bottom panels correspond to the solutions \eqref{eq:1D 1} and \eqref{eq:1D 2}, respectively.
As shown in the left panels, the RINN solution closely follow the exact one across the domain. The absolute error plots on the right support this result, showing that the majority of pointwise errors are confined to below $\mathcal{O}(10^{-8})$ in the first case and within $\mathcal{O}(10^{-6})$ in the second. 

\begin{example}  
    Consider 2D Poisson equations defined over the square domain $\Omega = [-1,1]^2$ with the following two analytical solutions:
    \begin{subequations}\label{eq:2D_all}
        \begin{align}
          u(\boldsymbol{x})&=\sin(6\pi x_1)\sin(6\pi x_2), \label{eq:2D 1} \\
          u(\boldsymbol{x}) &= -\prod_{i=1}^{2}  \left[2\cos\left(\frac{3}{2}\pi x_i + \frac{2\pi}{5}\right) + \frac{3}{2}\cos\left(3\pi x_i - \frac{\pi}{5}\right)\right]. \label{eq:2D 2} 
        \end{align}
    \end{subequations}
    \label{ep:Poi2D}
\end{example}

We employ a fully connected neural network with architecture $N_{\text{nn}} = [2, 512, 1024, 1]$. The number of collocation points is $K_{\mathrm{res}} = 2048$ and $K_{\mathrm{bcs}} = 4096$. Network weights are initialized from a uniform distribution $\mathcal{U}(-1, 1)$.
Training is performed using the Adam optimizer with a learning rate $\eta = 10^{-3}$ for $E = 500$ epochs. The regularization coefficient in the loss function \eqref{eq:loss} is set to $\varepsilon = 0.01$.

\begin{table}[htbp]
    \caption{Comparison of PIELM and RINN on \Cref{ep:Poi2D}}
    \label{table:2D test}
    \begin{tabular*}{0.75\textwidth}{@{\extracolsep{\fill}}cccc}
        \toprule
        Analytical Solution   & Method    & {$E_{L^2}$}   & {$E_{L^1}$}     \\
        \midrule
        \multirow{2}{*}{\eqref{eq:2D 1}} 
                    & PIELM & $8.74\times10^{-1}$ & $2.60\times10^{-1}$ \\
                    & RINN  & $1.22\times10^{-3}$ & $3.71\times10^{-4}$ \\
        \midrule
        \multirow{2}{*}{\eqref{eq:2D 2}} 
                    & PIELM & $9.40\times10^{-3}$ & $1.34\times10^{-2}$ \\
                    & RINN  & $1.06\times10^{-9}$ & $2.24\times10^{-9}$ \\
        \bottomrule
    \end{tabular*}
\end{table}

\Cref{table:2D test} presents a performance comparison between PIELM and RINN on two 2D Poisson equations with distinct exact solutions. In the first case \eqref{eq:2D 1}, PIELM yields a relative $L^2$ error around $\mathcal{O}(10^{-1})$ and a mean absolute error around $\mathcal{O}(10^{-1})$. In contrast, RINN significantly reduces these errors to $\mathcal{O}(10^{-3})$ and $\mathcal{O}(10^{-4})$, respectively.
For the second case \eqref{eq:2D 2}, RINN again achieves superior accuracy. The relative $L^2$ error decreases from $\mathcal{O}(10^{-2})$ with PIELM to below $\mathcal{O}(10^{-9})$ and the mean absolute error is similarly reduced from $\mathcal{O}(10^{-2})$ to $\mathcal{O}(10^{-9})$. 

\begin{figure}[htbp]
	\centering
	\includegraphics[width=0.90\linewidth]{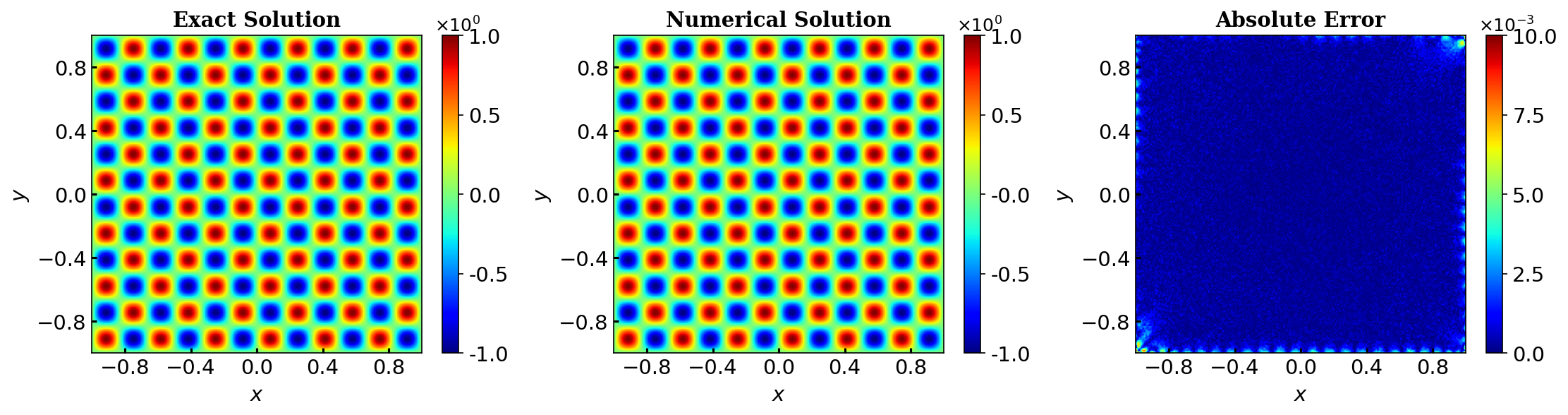}
	
	\includegraphics[width=0.90\linewidth]{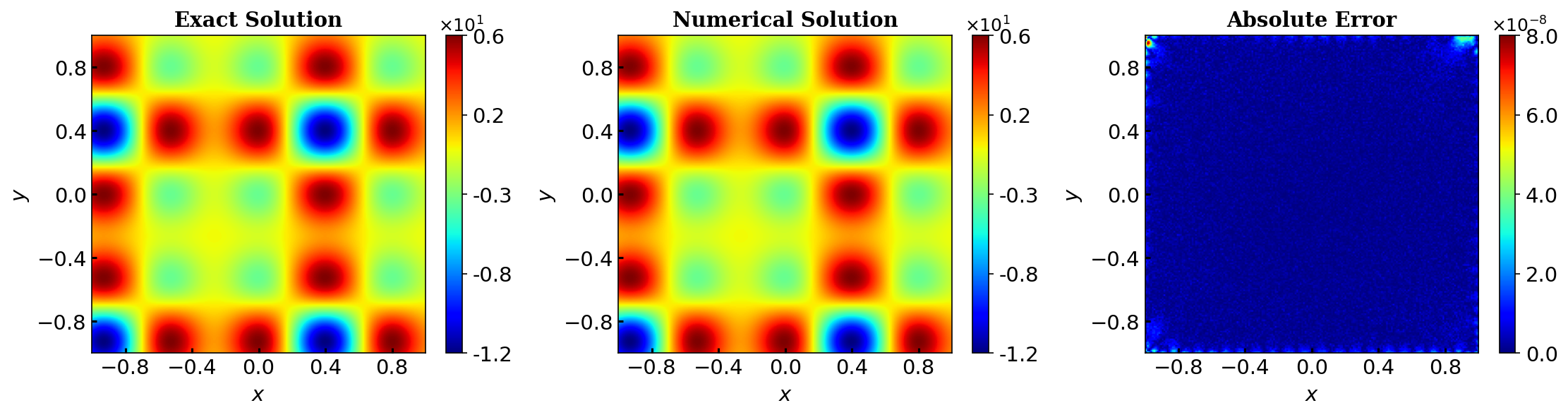}
    
    \caption{RINN results for \Cref{ep:Poi2D}. Each column shows: (left) analytical solution, (middle) predicted solution, (right) pointwise absolute error. Rows from top to bottom correspond to the high-frequency solution \eqref{eq:2D 1} and the multi-scale solution \eqref{eq:2D 2}.}
	\label{fig 2D}
\end{figure}

\Cref{fig 2D} shows the numerical results obtained by RINN for two 2D Poisson equations with distinct exact solutions. From left to right, each row displays the analytical solution, the RINN solution, and the corresponding pointwise absolute error. 
The error distributions in the right column indicate that RINN is capable of approximating both solutions with high accuracy. For \eqref{eq:2D 1}, the maximum pointwise error reaches $1 \times 10^{-2}$, while most errors are below  $\mathcal{O}(10^{-3})$. In the case of \eqref{eq:2D 2}, the error is significantly lower, with a maximum value around $8 \times 10^{-8}$ and the majority of values close to $\mathcal{O}(10^{-9})$.

\subsubsection{Evolution equations}
We apply the RINN to a series of time-dependent partial differential equations, including the advection equation, the heat equation, and the wave equation. Through systematic numerical experiments on these problems, we evaluated the numerical accuracy of this method and demonstrated its effectiveness and robustness in solving time-dependent partial differential equations.

To solve these problems using RINN, we adopt a fully connected neural network with architecture $N_{\text{nn}} = [2, 512, 1024, 1]$. All network weights are initialized from a uniform distribution $\mathcal{U}(-1, 1)$. The loss function defined in \eqref{eq:loss} is evaluated at residual points and optimized via the Adam optimizer with a constant learning rate of $\eta = 10^{-3}$ over $E = 500$ epochs. A regularization coefficient $\varepsilon = 0.01$ is employed during training.
Furthermore, the number of collocation points is $K_{\mathrm{res}} = 2048$, $K_{\mathrm{bcs}} = 2048$, and $K_{\mathrm{ics}} = 1024$.

\begin{example}
    (1D advection equations) Consider the following 1D advection equation defined on $\Omega=[-1,1], I=[0,1]$, 
    \begin{equation}
    \left\{
    	\begin{aligned}
    		u_t(x,t)&= -c u_x(x,t) ,&&\text{in } \Omega\times I, \\
    		u(-1,t)&=u(1,t),&&\text{in } I, \\
    		u(x,0)&=h(x),&&\text{on } \Omega\times \{0\}. \\
    	\end{aligned}
    \right.
    	\label{eq:Advection equations}
    \end{equation}
    Let the wave speed $c=0.4$ and two initial functions are defined as follows: 
    \begin{subequations}\label{eq:advection_all}
    \begin{align}
        h(x) &= e^{(0.5\sin(2\pi x))}-1,  \label{eq:Ad 1} \\
        h(x) &= \sin(2\pi x)\cos(3\pi x). \label{eq:Ad 2} 
    \end{align}
    \end{subequations}
    \label{ep:Adv1D}
\end{example}

\begin{table}[htbp]
    \caption{Comparison of PIELM and RINN on \Cref{ep:Adv1D}}
    \label{table:advection}
    \begin{tabular*}{0.75\textwidth}{@{\extracolsep{\fill}}cccc}
        \toprule
        Initial Condition    & Method      & {$E_{L^2}$}    & {$E_{L^1}$}     \\
        \midrule
        \multirow{2}{*}{\eqref{eq:Ad 1}} 
                    & PIELM & $8.28\times10^{-5}$ & $1.38\times10^{-5}$ \\
                    & RINN  & $1.85\times10^{-6}$ & $3.31\times10^{-7}$ \\
        \midrule
        \multirow{2}{*}{\eqref{eq:Ad 2}} 
                    & PIELM & $6.74\times10^{-4}$ & $1.59\times10^{-4}$ \\
                    & RINN  & $1.85\times10^{-8}$ & $6.03\times10^{-9}$ \\
        \bottomrule
    \end{tabular*}
\end{table}

\Cref{table:advection} presents a comparison between PIELM and RINN on two 1D advection problems. In the case of the smooth initial condition \eqref{eq:Ad 1}, RINN achieves a relative $L^2$ error of approximately $\mathcal{O}(10^{-6})$ and a mean absolute error of approximately $\mathcal{O}(10^{-7})$. For the second case \eqref{eq:Ad 2} which involves higher-frequency components, the advantage of RINN becomes more pronounced: the relative $L^2$ error reduces from $\mathcal{O}(10^{-4})$ to $\mathcal{O}(10^{-8})$ and the mean absolute error drops from $\mathcal{O}(10^{-4})$ to below $\mathcal{O}(10^{-8})$. 

\begin{figure}[htbp]
	\centering
	\includegraphics[width=0.90\linewidth]{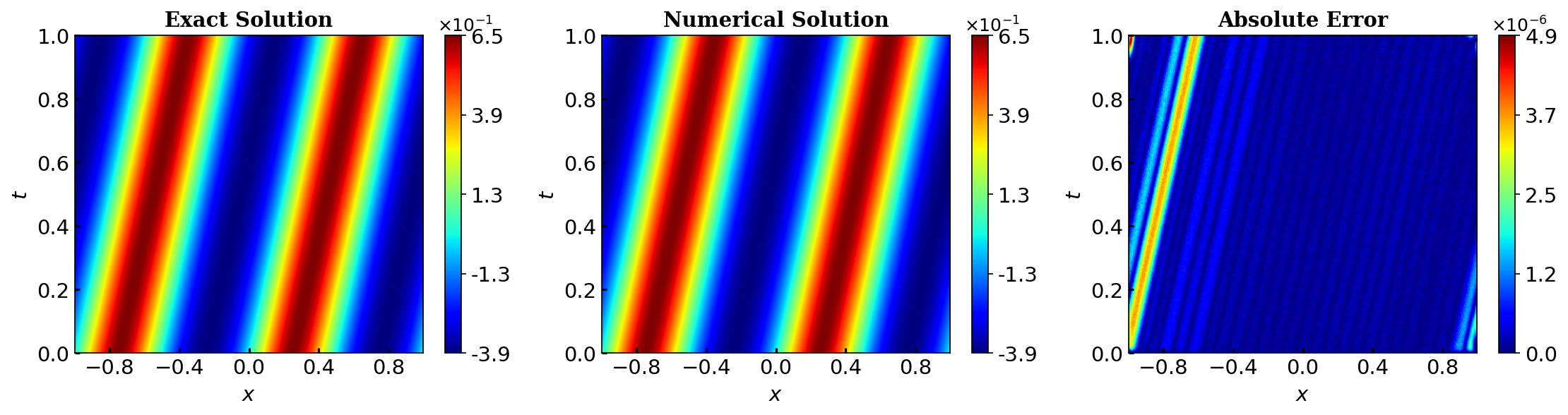}
	
	\includegraphics[width=0.90\linewidth]{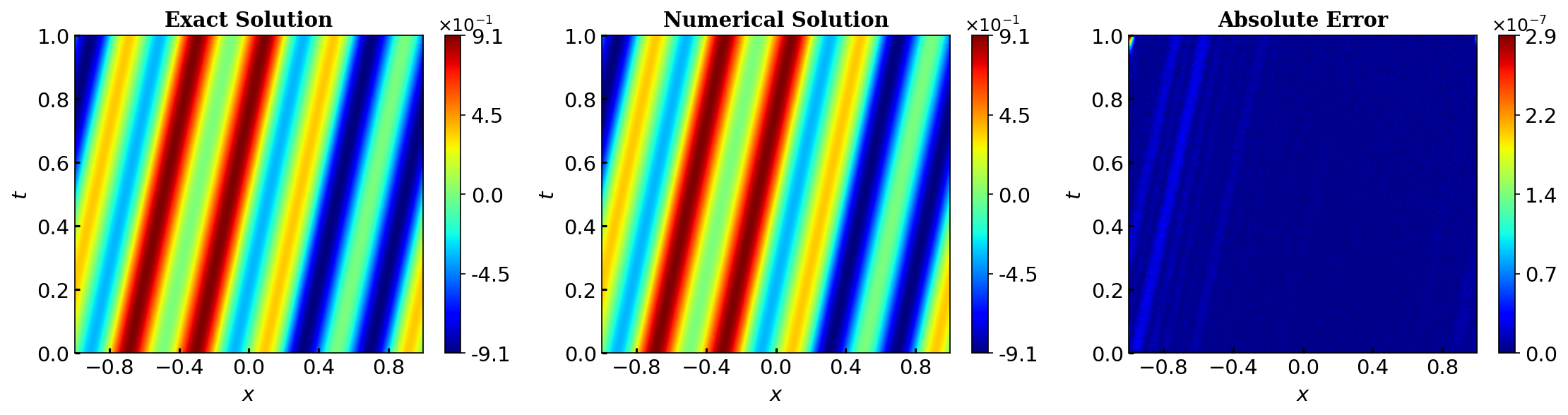}
	
	\caption{RINN results for \Cref{ep:Adv1D}. Each column shows: (left) analytical solution, (middle) predicted solution, (right) pointwise absolute error. Rows from top to bottom correspond to initial condition \eqref{eq:Ad 1} and \eqref{eq:Ad 2}.}
    \label{fig:advection}
\end{figure}

\Cref{fig:advection} shows the numerical results obtained by RINN for the 1D advection equation under the two initial conditions. In each row, the left column presents the analytical solution, the middle column displays the RINN solution, and the right column plots the pointwise absolute errors.
For the case of \eqref{eq:Ad 1} (top panel), the absolute error is primarily distributed along the transport direction and remains around $\mathcal{O}(10^{-6})$, with a maximum value of approximately $4.9 \times 10^{-6}$. In the second case \eqref{eq:Ad 2} (bottom panel), the error is further reduced to $\mathcal{O}(10^{-8})$, reaching a peak of about $2.9 \times 10^{-8}$. 

\begin{example}
	(1D heat equations) Consider the following 1D heat equation defined on $\Omega=[-1,1], I=[0,1]$, 
    \begin{equation}
    \left\{
    	\begin{aligned}
    		u_t(x,t) - u_{xx}(x,t)&=f(x,t) ,&& (x,t)\in \Omega\times I, \\
    		u(x,t)&=g(t),            && (x,t)\in \partial \Omega\times I, \\
    		u(x,t)&=h(x),            && (x,t)\in \Omega\times \{0\}. \\
    	\end{aligned}
    \right.
    	\label{eq:Heat equations}
    \end{equation}
    The analytical solution is defined as follows, 
	\begin{equation}
		u(x,t) = e^{-t}\sin(k\pi x),\qquad k=2, 6. \label{eq:time 1D}
	\end{equation}
    \label{ep:Heat1D}
\end{example} 

\begin{table}[htbp]
    \caption{Comparison of PIELM and RINN on \Cref{ep:Heat1D}}
    \label{table:heat}
    \begin{tabular*}{0.75\textwidth}{@{\extracolsep{\fill}}cccc}
        \toprule
        $k$       & Method                    & {$E_{L^2}$}        & {$E_{L^1}$}  \\
        \midrule
        \multirow{2}{*}{$2$} 
                    & PIELM & $2.91\times10^{-6}$ & $9.80\times10^{-7}$  \\
                    & RINN  & $6.49\times10^{-13}$ & $2.36\times10^{-13}$  \\
        \midrule
        \multirow{2}{*}{$6$} 
                    & PIELM & $1.27\times10^{-3}$ & $2.40\times10^{-4}$  \\
                    & RINN  & $2.48\times10^{-7}$ & $9.00\times10^{-8}$  \\
        \bottomrule
    \end{tabular*}
\end{table}

\Cref{table:heat} reports the numerical results for the 1D heat equation with the analytical solution given by \eqref{eq:time 1D}. Consider separately the cases where $k=2$ and where $k=6$.
When $k=2$, the relative $L^2$ error and mean absolute error of PIELM are both around $\mathcal{O}(10^{-6})$, while RINN reduces these to $\mathcal{O}(10^{-13})$. For $k=6$ which introduces more spatial variation, the relative $L^2$ error of PIELM is approximately $\mathcal{O}(10^{-3})$ and the mean absolute error is around $10^{-4}$. In contrast, RINN improves both metrics by roughly six orders of magnitude, achieving errors around $\mathcal{O}(10^{-7})$ and $\mathcal{O}(10^{-8})$, respectively.

\begin{figure}[htbp]
	\centering
	\includegraphics[width=0.90\linewidth]{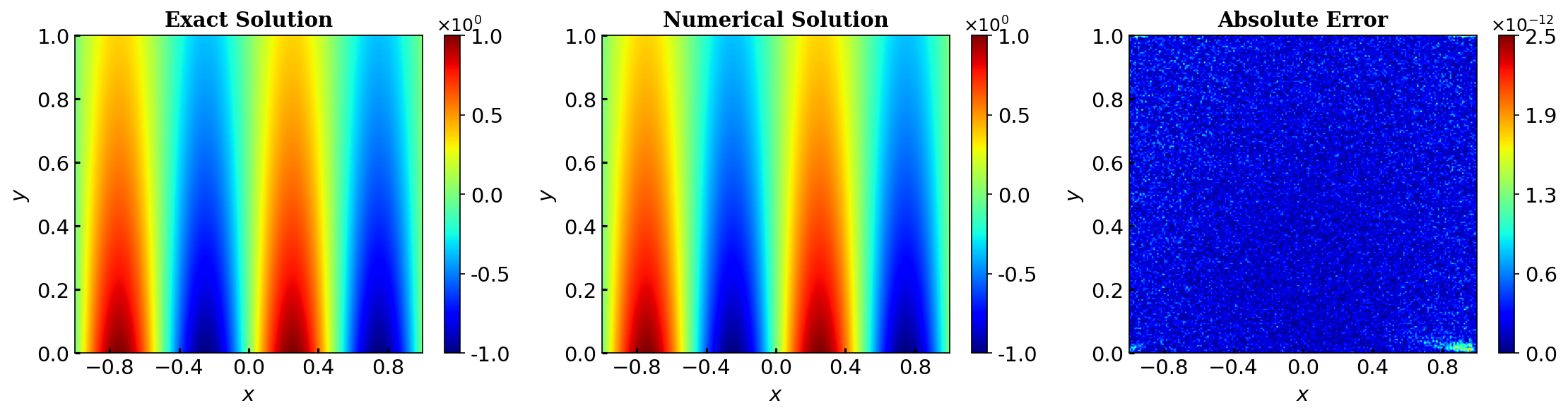}
	
	\includegraphics[width=0.90\linewidth]{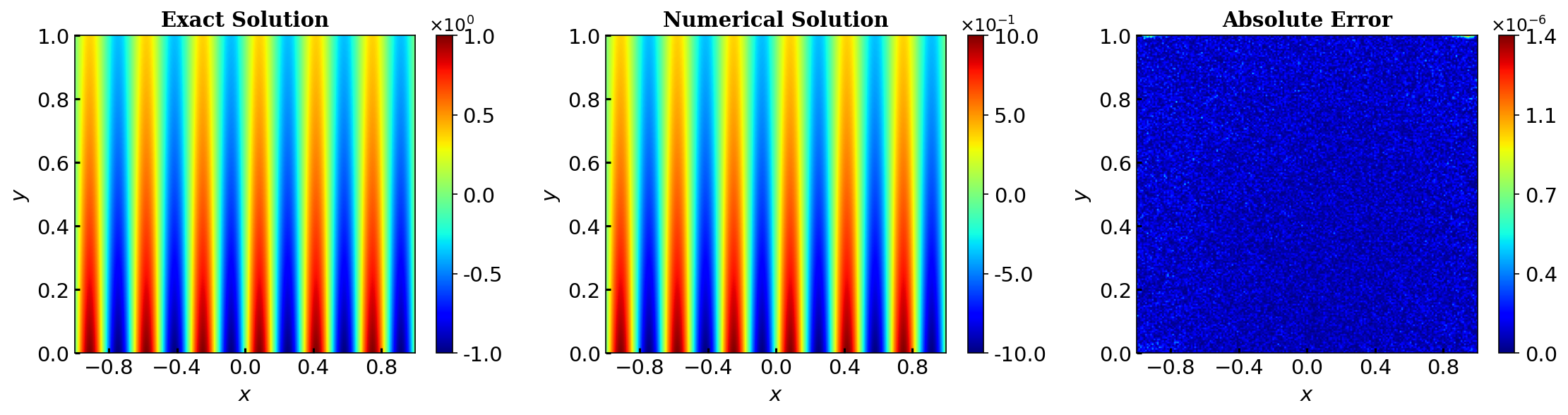}
	
	\caption{RINN results for \Cref{ep:Heat1D}. Each column shows: (left) analytical solution, (middle) predicted solution, (right) pointwise absolute error. Rows from top to bottom correspond to \eqref{eq:time 1D} with  $k=2$ and $k=6$.}
    \label{fig time1D}
\end{figure}

\Cref{fig time1D} presents the numerical results obtained by RINN for the 1D heat equation \eqref{eq:Heat equations} with the analytical solution given in \eqref{eq:time 1D}. Two wavenumber values are considered: $k=2$ (top row) and $k=6$ (bottom row). In each case, the left column shows the analytical solution, the middle column displays the RINN solution, and the right column illustrates the pointwise absolute error.
For $k=2$, the error remains consistently around $\mathcal{O}(10^{-12})$, with a maximum value of approximately $2.5 \times 10^{-12}$. The predicted solution aligns closely with the exact one over the entire domain.
For $k=6$, the solution contains more oscillatory features. In this case, the error increases to $\mathcal{O}(10^{-6})$, with the pointwise error attaining a maximum of $1.4 \times 10^{-6}$. While the solution exhibits finer spatial variations, the overall accuracy remains stable.

\begin{example}
    (1D wave equations) Consider the following 1D wave equation defined on $\Omega=[0,1], I=[0,1]$, 
    \begin{equation}
    \left\{
    	\begin{aligned}
    		u_{tt}(x,t)&= c^2 u_{xx}(x,t) ,&&   (x,t)\in \Omega\times I, \\
    		u(-1,t)&=u(1,t)=0,             &&   t\in I, \\
    		u(x,0)&=h(x),                  &&   x\in \Omega, \\
            u_t(x,0)&=0,                   &&   x\in \Omega. \\
    	\end{aligned}
    \right.
    	\label{eq:Wave equations}
    \end{equation}
    where $c=2$ is the wave speed and two initial functions are defined as follows,
    \begin{subequations}\label{eq:wave_all}
    \begin{align}
        h(x) &= \sin(\pi x) + \sin(2\pi x) + \sin(3\pi x),  \label{eq:Wave 1} \\
        h(x) &= \sin(2\pi x) + \sin(4\pi x). \label{eq:Wave 2} 
    \end{align}
    \end{subequations}
    \label{ep:Wave1D}
\end{example}

\begin{table}[htbp]
    \caption{Comparison of PIELM and RINN on \Cref{ep:Wave1D}}
    \label{table:wave}
    \begin{tabular*}{0.75\textwidth}{@{\extracolsep{\fill}}cccc}
        \toprule
        Initial Condition      & Method       & {$E_{L^2}$}     & {$E_{L^1}$}  \\
        \midrule
        \multirow{2}{*}{\eqref{eq:Wave 1}} 
                    & PIELM & $2.22\times10^{-4}$ & $1.55\times10^{-4}$ \\
                    & RINN  & $3.01\times10^{-8}$ & $2.07\times10^{-8}$  \\
        \midrule
        \multirow{2}{*}{\eqref{eq:Wave 2}} 
                    & PIELM & $8.20\times10^{-3}$ & $4.56\times10^{-4}$ \\
                    & RINN  & $1.03\times10^{-5}$ & $5.75\times10^{-6}$  \\
        \bottomrule
    \end{tabular*}
\end{table}

\Cref{table:wave} reports the numerical results for the 1D wave equation with two different initial conditions. Both cases are computed with wave speed $c=2$ and zero initial velocity. For \eqref{eq:Wave 1}, compared with the error of PIELM, the errors of RINN are reduced from the order of $\mathcal{O}(10^{-4})$ to the order of $\mathcal{O}(10^{-8})$. For \eqref{eq:Wave 2}, the relative $L^2$ error from PIELM is approximately $\mathcal{O}(10^{-2})$, while RINN lowers it to $\mathcal{O}(10^{-5})$. The mean absolute error also decreases from $\mathcal{O}(10^{-4})$ to $\mathcal{O}(10^{-6})$. 

\begin{figure}[htbp]
	\centering
	\includegraphics[width=0.90\linewidth]{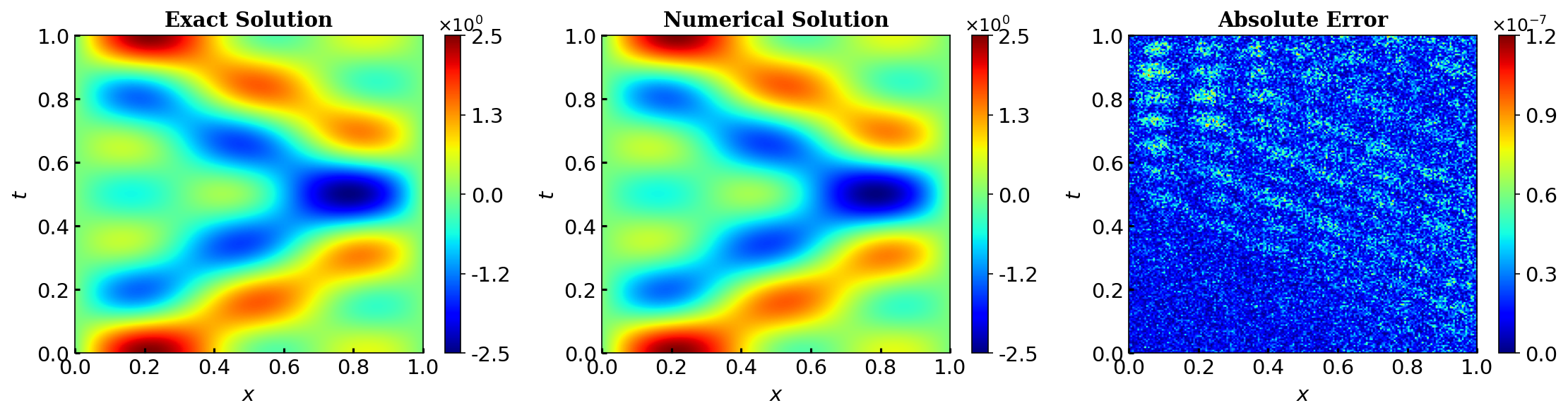}
	
	\includegraphics[width=0.90\linewidth]{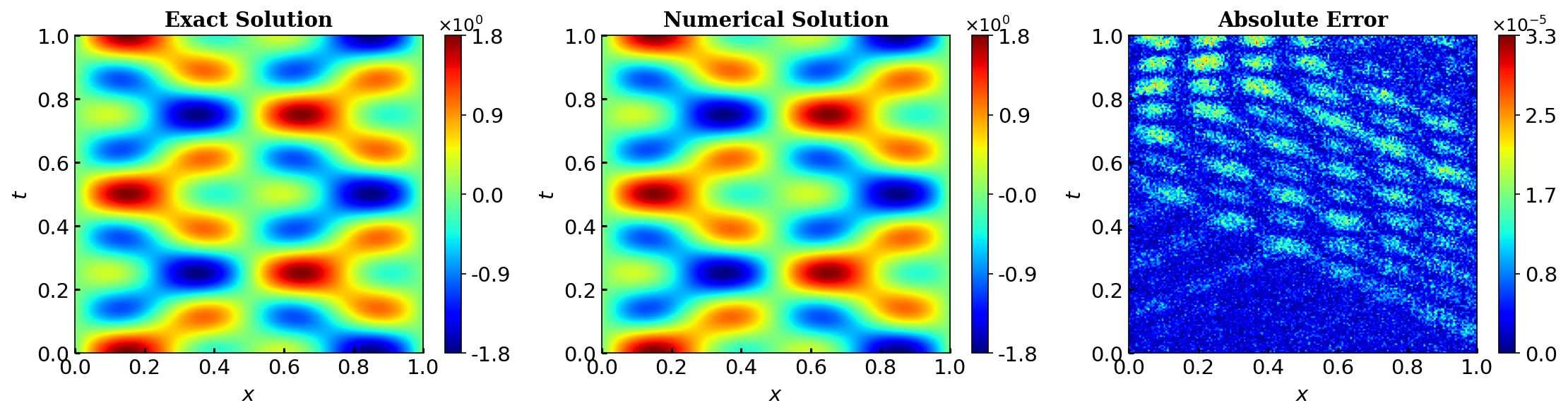}
	
	\caption{RINN results \Cref{ep:Wave1D}. Each column shows: (left) analytical solution, (middle) predicted solution, (right) pointwise absolute error. Rows from top to bottom correspond to initial condition \eqref{eq:Wave 1} and \eqref{eq:Wave 2}.}
    \label{fig:wave}
\end{figure}

\Cref{fig:wave} shows the results of RINN on the 1D wave equation with initial conditions. In each row, the left column presents the analytical solution, the center column displays the RINN prediction, and the right column plots the pointwise absolute error.
For the initial condition in \eqref{eq:Wave 1}, the absolute error remains $\mathcal{O}(10^{-7})$, with a maximum around $1.2 \times 10^{-7}$. the error distributes smoothly across the domain.
For the initial condition in \eqref{eq:Wave 2}, the error increases to $\mathcal{O}(10^{-5})$, with the pointwise error attaining a maximum of $1.4 \times 10^{-6}$. The error distribution follows the oscillatory structure of the solution.

\section{Rank Inspired Neural Network with Early Stopping}
\label{section:3}
    In the previous sections, we presented the covariance-driven orthogonalization strategy and demonstrated the effectiveness and stability of RINN across a range of PDE problems. To further enhance the final solution accuracy and mitigate potential degradation due to excessive training, this section begins by illustrating a possible increase in error when the degree of orthogonalization becomes overly high, using a simple one-dimensional Poisson equation as an example. Following this, we introduce an early stopping criterion based on PDE residuals, along with the detailed algorithmic implementation.

    \subsection{Early Stop Strategy Based On Residual Loss}

    In the following, we examine the issues caused by overtraining during the orthogonalization stage through a representative example and then introduce an early stopping strategy based on the residual loss.
    
    Consider 1D Poisson equation defined on $\Omega=[-1,1]$:
    $$
    \left\{
    \begin{aligned}
        -u_{xx}&=f, \quad &&x\in \Omega, \\
        u&=g,\quad &&x\in \partial \Omega.
    \end{aligned}
    \right.
    $$
    We take $\sin(2\pi x)\cos(3\pi x)$ as the analytical solution of the problem and employ a neural network with structure $N_{\text{nn}} = [1,128,1]$. The network weights are initialized using a uniform distribution $\mathcal{U}(-20,20)$. To construct the linear system, we use $K_{\mathrm{res}} = 1024$ residual points and $K_{\mathrm{bcs}} = 2$ boundary condition points. The RINN solution is then evaluated on a set of points sampled at equal intervals over the domain $[-1,1]$.

    \begin{figure}[htbp]
        \centering
        \includegraphics[width=0.6\linewidth]{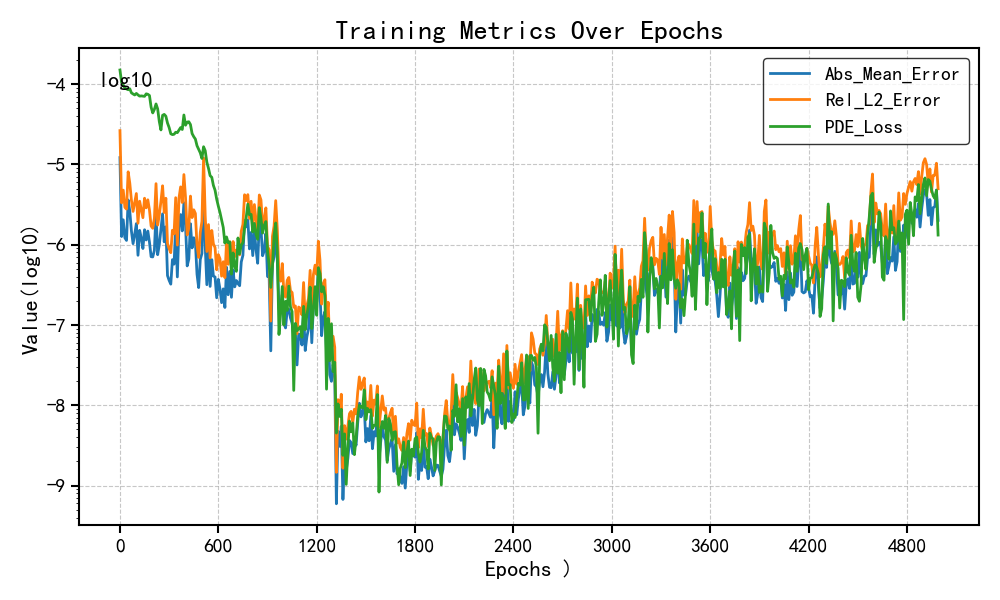}
        \caption{Evolution of the relative $L^2$ error and mean absolute error on evaluation points, alongside the PDE residual loss on training points.}
        \label{fig:Error_loss_epoch_1D}
    \end{figure}
    
    The PDE residual loss $\mathcal{L}_{\mathrm{pde}}$ is defined by
    $$
    \mathcal{L}_{\mathrm{pde}} = \mathcal{L}_{\mathrm{res}} + \mathcal{L}_{\mathrm{bcs}}
    $$
    where
    $$
    \mathcal{L}_{\mathrm{res}} = \sqrt{\frac{1}{K_{\mathrm{res}}} \sum_{i=1}^{K_{\mathrm{res}}} \bigl|u_{xx}(\mathbf{X}_{\mathrm{res}}^i) + f(\mathbf{X}_{\mathrm{res}}^i)\bigr|^2}, \quad
        \mathcal{L}_{\mathrm{bcs}} = \sqrt{\frac{1}{K_{\mathrm{bcs}}} \sum_{i=1}^{K_{\mathrm{bcs}}} \bigl|u(\mathbf{X}_{\mathrm{bcs}}^i) - g(\mathbf{X}_{\mathrm{bcs}}^i)\bigr|^2}.
    $$
    
    \Cref{fig:Error_loss_epoch_1D} shows the relative $L^2$ error, mean absolute error, and PDE residual loss plotted against training epochs. All three metrics initially decrease but then exhibit an upward trend after a certain point which may be related to excessive decorrelation of the basis functions or parameter overfitting.
    The PDE residual loss $\mathcal{L}_{\mathrm{pde}}$ is directly accessible during training and he variation of $\mathcal{L}_{\mathrm{pde}}$ closely follows that of the evaluation errors. Consequently, it can serve as an effective early stopping criterion. Terminating the orthogonalization training when $\mathcal{L}
    _{\mathrm{pde}}$ attains its minimum helps to prevent deterioration of the final solution accuracy.
    More generally, for PDE problems \eqref {eq1}, sampling points set is defined as $\mathcal{X} = \{\mathbf{X}_{\mathrm{res}},\,\mathbf{X}_{\mathrm{bcs}},\,\mathbf{X}_{\mathrm{ics}}\}$. The corresponding PDE residual losses are defined as follows:
    \begin{equation}
        \begin{aligned}
        \mathcal{L}_{\mathrm{pde}} &= \mathcal{L}_{\mathrm{res}}+\mathcal{L}_{\mathrm{bcs}}+\mathcal{L}_{\mathrm{ics}},
        \end{aligned}
        \label{eq:pde loss}
    \end{equation}
    where
    $$
    \begin{aligned}
        \mathcal{L}_{\mathrm{res}} &=\sqrt{\frac{1}{K_{\mathrm{res}}} \sum_{i=1}^{K_{\mathrm{res}}} \left| \mathscr{A}u(\mathbf{X}_{\mathrm{res}}^i) - f(\mathbf{X}_{\mathrm{res}}^i) \right|^2}, \\
        \mathcal{L}_{\mathrm{bcs}} &=\sqrt{\frac{1}{K_{\mathrm{bcs}}} \sum_{i=1}^{K_{\mathrm{bcs}}} \left| \mathscr{B}u(\mathbf{X}_{\mathrm{bcs}}^i) - g(\mathbf{X}_{\mathrm{bcs}}^i) \right|^2}, \\
        \mathcal{L}_{\mathrm{ics}} &=\sqrt{\frac{1}{K_{\mathrm{ics}}} \sum_{i=1}^{K_{\mathrm{ics}}} \left| u(\mathbf{X}_{\mathrm{ics}}^i) - h(\mathbf{X}_{\mathrm{ics}}^i) \right|^2}.
    \end{aligned}
    $$
    Based on the residual loss defined in \eqref{eq:pde loss}, we propose an early stopping strategy at the algorithmic level, as described in \Cref{alg:RINN_es}. This algorithm retains the original two-stage framework of orthogonalization and least squares, while adding residual evaluation for early stopping.
    
\begin{algorithm}[]
  \caption{Rank‑Inspired Neural Network with Early Stopping (RINN-es)}
  \label{alg:RINN_es}
  \KwIn{Collocation points $\mathcal{X} = \{\mathbf{X}_{\mathrm{res}}, \mathbf{X}_{\mathrm{bcs}}, \mathbf{X}_{\mathrm{ics}}\}$, 
        maximum epochs $E$, learning rate $\eta$, regularization coefficient $\varepsilon$, max patience $P$.}
  \KwOut{Optimized parameters $(\boldsymbol{\theta}, \boldsymbol{\beta})$.}
 
  Initialize network parameters $\boldsymbol{\theta}$ randomly\;
  Set tolerance $\texttt{tol}= +\infty$, 
      no improve counter $\texttt{no\_improve}= 0$\;
  Initialize Adam optimizer with learning rate $\eta$\;

  \For{$epoch := 1$ \KwTo $E$}{
    \textbf{Stage 1}: \texttt{Covariance‑Driven Basis Orthogonalization} \;
    \, Compute hidden layer output matrix $\boldsymbol{\Phi} = \boldsymbol{\Phi}_{\boldsymbol{\theta}}(\mathcal{X})$\;
    \, Compute covariance $\mathbf{C} = \tfrac{1}{K-1}\,\boldsymbol{\Phi}^\top \boldsymbol{\Phi}$\;
    \, Compute$\mathcal{L}_{\mathrm{total}} = \varepsilon\,\mathcal{L}_{\mathrm{diag}} + \mathcal{L}_{\mathrm{ortho}}$\;
    \, Update $\boldsymbol{\theta} = \boldsymbol{\theta} - \eta\,\nabla_{\boldsymbol{\theta}}\,\mathcal{L}_{\mathrm{total}}$\;

    \textbf{Stage 2}: \texttt{Least‑Squares Solve}\;
    \, Assemble linear system $\boldsymbol{H}\,\boldsymbol{\beta} = \boldsymbol{S}$\;
    \, Solve $\boldsymbol{\beta} = \boldsymbol{H}^\dagger\,\boldsymbol{S}$\;

    \textbf{Stage 3}: \texttt{PDE Residual and Early Stopping}\;
    \, Compute PDE residual loss $\mathcal{L}_{\mathrm{pde}}$\;
    \, \uIf{$\mathcal{L}_{\mathrm{pde}} < \texttt{tol}$}{
      $\texttt{tol} = \mathcal{L}_{\mathrm{pde}}$;\quad 
      save $\boldsymbol{\theta}^* = \boldsymbol{\theta}$, 
      $\boldsymbol{\beta}^* = \boldsymbol{\beta}$;\quad 
      $\texttt{no\_improve} = 0$\;
    }
    \, \Else{
      $\texttt{no\_improve} = \texttt{no\_improve} + 1$\;
      \If{$\texttt{no\_improve} \ge P$}{
        break 
      }
    }
  }

  \If{$\theta^*$ exists}{
    Restore $\boldsymbol{\theta} = \boldsymbol{\theta}^*$, 
    $\boldsymbol{\beta} = \boldsymbol{\beta}^*$\;
  }

  \KwRet{$(\boldsymbol{\theta}, \boldsymbol{\beta})$}.
\end{algorithm}

\subsection{Numerical Examples}

    In this subsection, we revisit the 2D Poisson problem \eqref{eq:Elliptic Equations} with the high-frequency solution \eqref{eq:2D 1} to validate the effectiveness of \Cref{alg:RINN_es} through an example. The network structure and the number of collocation points were all consistent with those in \Cref{ep:Poi2D}. 

    To evaluate the robustness of PIELM, RINN and RINN-es in terms of weight initialization, we respectively took $\mathcal{U}(-0.5, 0.5)$, $\mathcal{U}(-1, 1)$ and $\mathcal{U}(-2, 2)$.
    Both RINN and RINN-es employ the Adam optimizer with a learning rate of $\eta = 10^{-3}$ and both have a regularization coefficient of $\varepsilon = 0.01$. The RINN model is trained for a fixed number of $E = 500$ epochs, while RINN-es allows up to $E = 1000$ epochs, implementing an early stopping strategy with a maximum patience of $P = 250$ epochs.

    \begin{table}[htbp]
		\centering
        \caption{Error Comparison and Early Stopping Performance on 2D Poisson Equation under Different Initializations ($E$: maximum iterations, $E^*$: optimal iterations).}
        \label{table:es_Poi2D}
		\begin{tabular*}{0.75\textwidth}{@{\extracolsep{\fill}}ccccc}
			\toprule
			initialization & Method & $E$/$E^*$ & $E_{L^2}$ & $E_{L^1}$  \\
			\midrule
			\multirow{3}{*}{$\mathcal{U}(-0.5,0.5)$} 
    			& PIELM   &-/-    & $7.80\times10^{-2}$  & $3.17\times10^{-2}$ \\
    			& RINN    &500/-  & $1.16\times10^{+1}$  & $4.60\times10^{+0}$ \\
    			& RINN-es &1000/6 & $8.41\times10^{-3}$  & $3.34\times10^{-3}$ \\
			\midrule
			\multirow{3}{*}{$\mathcal{U}(-1,1)$}     
    			& PIELM   &-/-      & $8.74\times10^{-1}$  & $2.60\times10^{-1}$  \\
    			& RINN    &500/-    & $1.22\times10^{-3}$  & $3.71\times10^{-4}$  \\
    			& RINN-es &1000/672 & $9.57\times10^{-4}$  & $2.95\times10^{-4}$  \\
			\midrule
			\multirow{3}{*}{$\mathcal{U}(-2,2)$}     
    			& PIELM   &-/-      & $8.36\times10^{+0}$  & $3.43\times10^{+0}$  \\
    			& RINN    &500/-    & $3.46\times10^{+0}$  & $9.04\times10^{-1}$  \\
    			& RINN-es &1000/999 & $3.04\times10^{-2}$  & $6.50\times10^{-3}$  \\
			\bottomrule
		\end{tabular*}
        
	\end{table}

    \Cref{table:es_Poi2D} reports relative \(L^2\) and mean absolute errors for PIELM, RINN and RINN-es under three initialization ranges. Using \(\mathcal{U}(-0.5,0.5)\), RINN diverges with \(\mathcal{O}(10^{+1})\) errors while RINN-es converges by epoch~6 with \(\mathcal{O}(10^{-3})\) errors. Under \(\mathcal{U}(-1,1)\), RINN achieves \(\mathcal{O}(10^{-3})\) and RINN-es reduces the error to \(\mathcal{O}(10^{-4})\) by epoch~672. With \(\mathcal{U}(-2,2)\), PIELM and RINN degrade to \(\mathcal{O}(10^{0})\) and \(\mathcal{O}(10^{+1})\), respectively, whereas RINN-es maintains \(\mathcal{O}(10^{-2})\) at epoch~999. The optimal epochs (6, 672 and 999) show that fixed training lengths are inadequate. These results confirm that PDE-residual-based early stopping dynamically adapts to different initializations.

    \begin{figure}[htbp]
    	\centering
    	\includegraphics[width=0.90\linewidth]{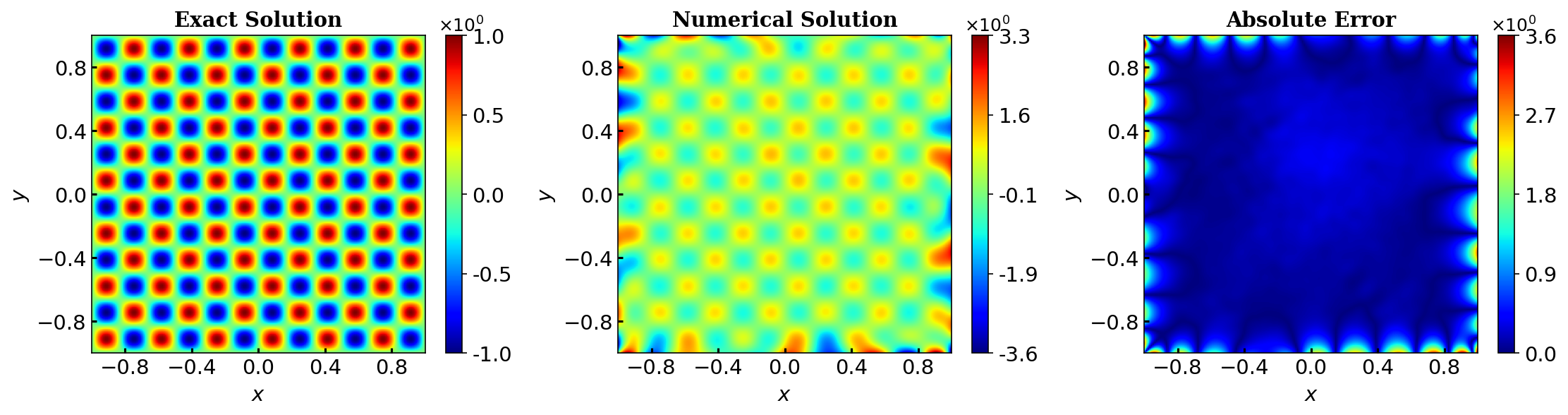}
    	
    	\includegraphics[width=0.90\linewidth]{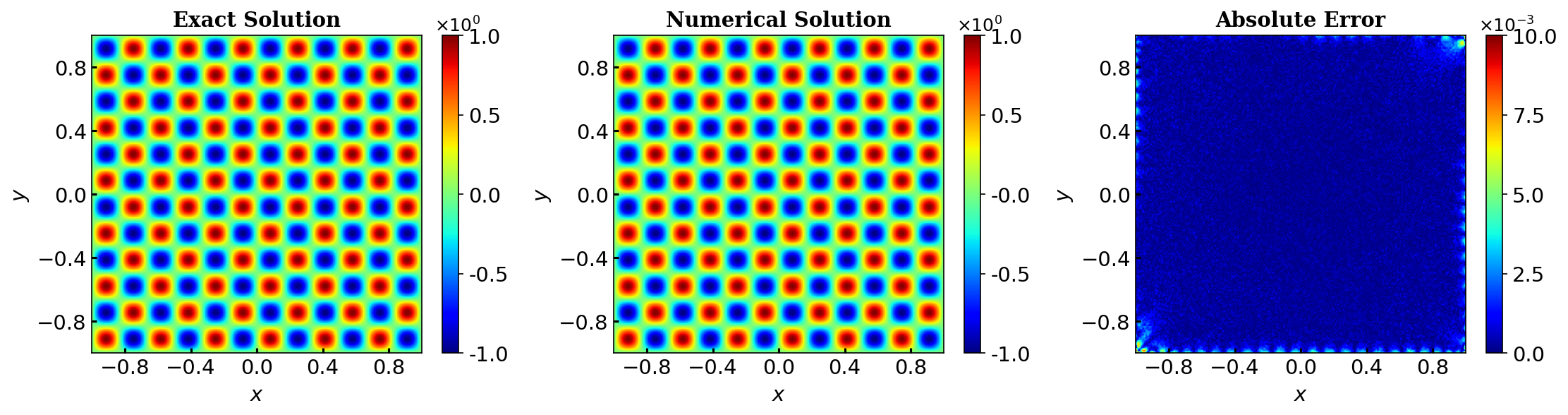}

        \includegraphics[width=0.90\linewidth]{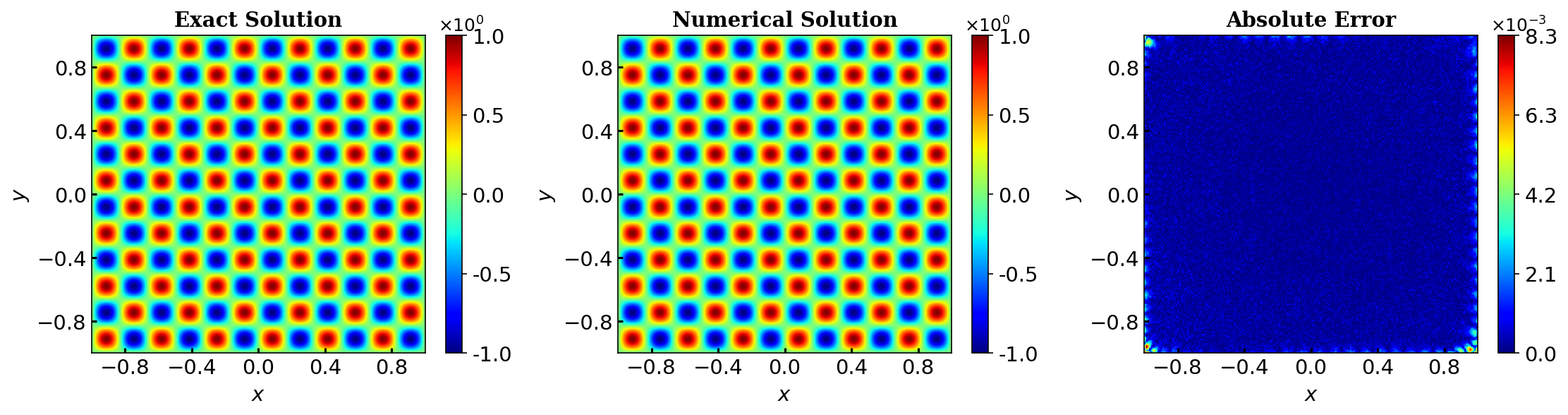}
    	
    	\caption{Numerical solutions and pointwise absolute error distributions for the 2D Poisson equation with initialization $\mathcal{U}(-1,1)$: PIELM (top panel), RINN (middle panel), and RINN-es (bottom panel).}
        \label{fig:es_Poi2D}
    \end{figure}

    \Cref{fig:es_Poi2D} shows the numerical solutions and pointwise absolute error distributions obtained by PIELM (top), RINN (middle), and RINN-es (bottom) for the 2D Poisson problem under initialization \(\mathcal{U}(-1,1)\). 
    In the top panel, PIELM produces large boundary errors, with absolute values reaching \(\mathcal{O}(10^0)\) near the domain edges. The interior solution also fails to capture the expected oscillatory structure and deviates significantly from the exact solution.
    The middle panel displays the RINN result. The solution is smoother and better represents the oscillatory pattern. However, boundary errors remain, with maximum values around \(\mathcal{O}(10^{-3})\).
    In the bottom panel, RINN-es shows further improvement due to early stopping. The predicted solution closely matches the reference throughout the domain. The maximum absolute error stays below \(\mathcal{O}(10^{-3})\), even near the boundaries.

\section{Conclusion}
\label{section:4}
    The paper introduces the Rank Inspired Neural Network (RINN) framework to improve the robustness of PIELM by addressing instability caused by strongly correlated hidden-layer basis functions. RINN employs a covariance-driven orthogonalization pretraining stage, separating the orthogonalization of neural basis functions from solving the linear equations. This creates a more stable two-stage training procedure. An early-stopping criterion based on PDE residual monitoring is integrated to prevent excessive decorrelation and maintain approximation fidelity.
    
    Numerical experiments on elliptic and evolution PDEs confirm RINN's effectiveness. It maintains high accuracy across diverse problem scales and physical regimes. Compared to conventional PIELM, RINN shows significantly improved initialization robustness.
    
    Future work will extend to systems of nonlinear PDEs, irregular domain geometries, and multiscale problems.

\section*{Acknowledgments}
This work was supported by NSFC Project (12431014) and Project of Scientiﬁc Research Fund of the Hunan Provincial Science and Technology Department (2024ZL5017).
	

\end{document}